\documentclass[final,11pt,onecolumn,a4]{IEEEtran}

\newcommand{\bP}{\ensuremath{\mathbf{P}}}
\newcommand{\bQ}{\ensuremath{\mathbf{Q}}}

\newcommand{\eye}{\ensuremath{\mathbf{I}}}
\newcommand{\bH}{\ensuremath{\mathbf{H}}}

\newcommand{\bx}{\ensuremath{\mathbf{x}}}
\newcommand{\by}{\ensuremath{\mathbf{y}}}

\newtheorem{definition}{Definition}[section]

\newcommand{\cov}{\ensuremath{\mathrm{cov}}}

\newcommand{\x}{{\bf x}}

\newcommand{\bxi}{{\bm{\xi}}}
\newcommand{\bnu}{{\bm{\nu}}}

\newcommand{\y}{{\bf y}}
\newcommand{\be}{{\bf b}}

\newcommand{\bomega}{{\bm{\omega}}}

\newcommand{\rott}{{\mathbf{r}}_{-\theta}}

\newcommand{\bj}{\bm{j}}

\usepackage{graphics,graphicx,bm,amsmath,amssymb}      
\usepackage{pst-all,amstext,amsfonts,epsf,epsfig,pifont}
\begin{document}
%
\title{Localisation of Geometric Anisotropy}
%
%
\author{Sofia~C.~Olhede
\thanks{S. Olhede is with Imperial College London,
SW7 2AZ, London, UK (s.olhede@imperial.ac.uk). Tel:
+44 (0) 20 7594 8568, Fax: +44 (0) 20 7594 8517.}}

\markboth{Statistics Section Technical Report TR-07-01}{Olhede:
Localisation of Geometric Anisotropy}
\maketitle

\begin{abstract}
The class of 2-D nonseparable geometrically anisotropic localisation operators
is defined, containing highly anisotropic nearly unidirectional localisation
operators,
as well as isotropic localisation operators. A continuum
of anisotropic operators between the extremes of near unidirectionality and
isotropy
are treated in a single
class. 
The
eigensystem of any given operator in this family
is determined, thus specifying geometrically anisotropic optimally concentrated
functions, and their degree of localisation.
\end{abstract}
\vspace{0.1in}

\section{Introduction \label{intro}}
\PARstart{T}his correspondence introduces a new class
of spatially anisotropic and nonseparable 2-D localisation operators, namely
the class of {\em geometrically
anisotropic}
localisation operators. The operation of `localisation' in this context refers to limiting
variation in a 2-D square integrable function to variation associated
with a given set of spatial locations and spatial frequencies, see also Daubechies \cite{Daub}[ch.~2]. By formally using a localisation
operator, the energy concentration
of a
given function to an anisotropic and nonseparable 4-D region of 2-D space
and 2-D spatial
frequency, may
be {\em exactly} quantified. Optimally localized functions can then be derived.

The important set of tools which motivates the need for deriving optimally localised
functions
is that based on the local representation of functions. A Fourier Transform
(FT) represents
a signal globally, {\em i.e.} it decomposes a signal in terms of modes present
over the full length of a signal. Given many signals exhibit transient features,
it is necessary to represent the local properties of the signal. Examples of local decompositions include the Wavelet
Transform (WT) and the windowed FT \cite{AMWA}.

To be able to form a local decomposition of a given function, a 
set of well-localised
decomposition functions
must be used, or the utility
of the decomposition vanishes. The locality of a given decomposition function
must be adjusted to the class of signal that will be analysed.
By
the construction procedure used for 
the WT \cite{AMWA}, the locality of any member of the family
of functions used for analysis, is in some sense `equivalent' to the locality
of the mother wavelet. Thus it is sufficient to determine appropriate mother
wavelets for a given problem.
The
eigensystem of a localisation operator corresponds
to a whole set of mutually orthogonal and optimally localized functions, where the
eigenvalues
of the functions, measure their degree of localisation to a given localisation
region
\cite{Metikas,Daub}. Thus eigenfunctions of localisation operators are suitable
mother wavelets.

By determining the eigensystem of the operators proposed in this correspondence, we obtain classes of 
mother wavelets suitable for analysis of non-stationary geometrically anisotropic
fields.
These facts justify the construction and study of the proposed operators.
1-D localisation operators have already
been
the focus of considerable study, see for example the references in 
\cite{Daub,DaubechiesPaul1988,Matz2006}.
Consecutive truncations form a possible method of constructing localisation
operators: such procedures
treat the space and spatial frequency variables inhomogeneously. This results in
an unequal compromise between concentration in space and in spatial frequency,
a clearly undesirable feature.

Also it is desirable to define operators localising functions to regions
that do not contain zero frequency,
if we wish to build for example families of wavelets.
In this case the relocation in frequency is not implemented by frequency
shifts.
Daubechies and Paul \cite{DaubechiesPaul1988} defined a set of homogeneous
localisation operators, that treated time and frequency variation on an equal
footing. These operators were used
to derive optimal 1-D mother wavelet functions, and the 1-D Morse
wavelets were thus obtained. The Morse wavelets
have been used to analyse geophysical, astrophysical, and medical time
series \cite{Olhede2,olhedemed}, and have been extended to radial wavelets
in 2-D \cite{Metikas}. The construction of localisation operators in both
cases started from forming a local decomposition of the observed function,
discussed in section \ref{local}.

We seek to develop operators for the study of geometric anisotropy, formally
defined in section \ref{geoani}. Our motivation for studying
this class is that it contains isotropic functions, and anisotropic functions
that are nearly unidirectional, as well as the
full continuum of structures in between isotropy and unidirectionality. Geometric
anisotropy is used for modelling covariance structures
in geostatistics (see for example Christakos \cite{Christakos}[p.~61]). Thus the optimally localised 
geometrically anisotropic mother wavelets will have
a natural application area, in the study of non-stationary geometrically
anisotropic random fields.

In general it is very
hard to define operators so that their eigensystem may be determined 
\cite{Daub}[p.~41], and
without a known eigensystem the operators lose most of their utility.  We shall therefore
define the
anisotropic operator with great care, in a series of steps, in section \ref{conop}.
We start by defining the appropriate `fiducial vector' from which to build the operator.
Geometric anisotropy is characterised by a transformation matrix, and thus a transformed
distance metric is defined. The fiducial vector is picked to be a function marginally well
localized in space in terms of the transformed distance metric, and also well-concentrated in transformed frequency. From
the fiducial vector a family of coherent states must subsequently be defined.
We construct the family of coherent states in section \ref{defstate}, using the transformation matrix of the geometric anisotropy. In general the family
of coherent states will not be constructed
by the usual 2-D operations of translating, scaling and rotating the fiducial vector \cite{Metikas}, but rather by a
set of {\em transformed} operations. 

We subsequently need to determine the localisation of
an arbitrary element of the family of coherent states, as this in conjunction
with the choice of region of integration, will determine the localisation
of the operator, see section \ref{locally}. This enables the definition of operators with given
localisation regions, see section \ref{defop}. We subsequently determine
the eigensystem of the proposed localisation operators, see section \ref{eigen}. We derive explicit forms for both eigenfunctions
and eigenvalues. We give specific
examples, for given transformation matrices, of functions in this class.

These developments combine to define new classes of `optimally' localised nonseparable anisotropic functions. The developments span previously derived optimally
localised functions, and form a cohesive framework for treating geometrically
anisotropic 2-D localisation.

\section{Local Decompositions \label{local}}
We shall construct the localisation operator starting from a local decomposition,
similar to the 2-D Continuous
WT (CWT), see for example Antoine {\em et al.} \cite{AMWA}. 
We define a family of `coherent states'
$v_{\bm{\xi}}(\x)$ from the `fiducial vector' \cite{Metikas}, a single 2-D function $v(\x)$. $v(\x)$ is
assumed to be marginally well-localised in space and spatial frequency. The index $\bm{\xi}=\left[a,\theta,\be\right]^T\in {\cal A}\subset {\mathbb{R}}^4$ regulates the localisation of the elements of the family
of coherent states. 
For example, the usual WT corresponds to taking:
$v_{\bm{\xi}}(\x)=
{\cal D}_a {\cal{R}}_{\theta} {\cal T}_{\be}\left\{
v\right\}(\x),
$
where for $a>0,$ ${\cal D}_a v(\x)=a^{-1}v(\x/a)$ is the dilation operator,
for $\theta\in\left[0,2\pi\right)$,
${\cal R}_{\theta}v(\x)=v\left(\rott \x\right)$ is the rotation operator, with ${\mathbf{r}}_{\theta}$ as the rotation matrix \cite{AMWA},
and for $\be \in {\mathbb{R}}^2,$ ${\cal T}_{\be}v(\x)=v(\x-\be)$
is the translation operator.

Let $v(\x)$ have a FT given by
$V(\bomega),$ in angular frequency $\bomega$. We do not refer to $v(\x)$ as a `mother wavelet', as we intend
to introduce a different set of operations to construct $v_{\bm{\xi}}(\x)$.
We define the `local coefficients' of $v,g\in L^2({\mathbb{R}}^2)$ by:
\begin{eqnarray}
w_{v}\left(\bm{\xi};g\right)&=&
\int \int_{{\mathbb{R}}^2} v_{\bm{\xi}}^*(\x)g(\x)\;d^2\x=\langle
v_{\bm{\xi}}(\x),g(\x)\rangle=\frac{1}{(2\pi)^2}\int \int_{{\mathbb{R}}^2}
V_{\bm{\xi}}^*(\bomega)G(\bomega)\;d^2\bomega
\label{2dcwt}.
\end{eqnarray}
Eqn. (\ref{2dcwt}) represents a projection of $g(\x)$ into `local contributions'
$w_{v}\left(\bm{\xi};g\right).$
The function $g(\x)$ can be reconstructed from the local coefficients,
if $\bm{\xi}$ and $v(\x)$ satisfy a suitable set of constraints. For example
if $v(\x)$ is an admissible mother wavelet, with admissibility constant $C_v$ and $\bm{\xi}$ corresponds to
the indexing denoting the translation,
scaling and rotation operation, then we may reconstruct $g(\x)$ by:
\begin{equation}
g(\x)=\frac{1}{C_{v}}\int_{\cal A} v_{\bm{\xi}}(\x)
w_{v}\left(\bm{\xi};g\right)\;dA_{\bm{\xi}},\quad dA_{\bm{\xi}}=
a^{-3}da\;d^2\be\;d\theta
\label{2dicwt}.
\end{equation}
The interpretability of eqns. (\ref{2dcwt}) and (\ref{2dicwt}) depends
on the locality of $v_{\bm{\xi}}(\x),$ in turn determined
from the locality of $v(\x)$ and the choice of operations denoted by $\bm{\xi}$.
Depending on the type of function $g(\x)$ that we are decomposing, different
families of coherent states are suitable to use for the analysis of this
function. To focus our interest on a special class of anisotropy, we shall
now introduce the class of geometric anisotropy.

\section{Constructing the Geometrically Anisotropic Operator 
\label{conop}}
\subsection{Geometrically Anisotropic Functions \label{geoani}}
\begin{definition}{A Geometrically Anisotropic Function}\\
A function $g_A(\x)$ is said to exhibit {\em geometric anisotropy}
if for a fixed non-negative $2 \times 2$ symmetric matrix $\bH$ and $g_r(\cdot)$
a 1-D function, it
takes the form:
$
g_A(\x)=g_r\left(\x^T \bH \x\right).
$
\end{definition}
An application of such functions is to model auto-covariance of
random fields having geometric anisotropy, see Christakos 
\cite{Christakos}[p.~61]. For such random fields the auto-covariance of the field
$h(\x)$
at spatial locations $\x\in {\mathbb{R}}^2$ and at $\y\in {\mathbb{R}}^2$ takes the form $b_h\left(\x,\y\right)=\cov\left\{h(\x),
h(\y)\right\}=b_r\left(\left(\x-\y\right)^T \bH \left(\x-\y\right)\right)$.
An example of such functions is an 
isotropic covariance, given by
$\bH=\eye,$ where $b_{\bH}(\x_1,\x_2)$ is specified as a function of the Cartesian distance
between $\x_1$ and $\x_2$, but also for certain observed phenomena it is natural to model
the auto-covariance between the field at two points as depending on a local
affine transformation of the two variables \cite{Christakos,Christakos+2000}. Note that the anisotropic extreme of this form
corresponds to having
$\bH=\eye_1(\varepsilon)=\left(\left[1\quad 0 \right];
\left[0\quad \varepsilon^2 \right]\right)$. The function becomes approximately unidirectional,
as $\varepsilon\rightarrow 0.$ 

We form the decomposition of $\bH$ as
$\bH=\bP_{\bH}^T\bP_{\bH}$, and refer to $\bP_{\bH}$ as the `transformation
matrix' of $\bH$.
Unfortunately this specification is not sufficiently constrained to
uniquely determine $\bP_{\bH}$, however as we in this correspondence only
seek to determine the optimal decomposition for a set $\bH$, any of the $\bP_{\bH}$ will work.
For $\bH=\eye,$ we have $\bP_{\bH}=\eye$ 
whilst for $\bH=\eye_1(\varepsilon),$
$\bP_{\bH}(\varepsilon)=\left(\left[1\quad 0 \right];
\left[0\quad \varepsilon \right]\right).$ Let ${\mathbf Q}_{\bH}={\mathbf
P}^{-1}_{\bH}$.
Define the `transformed spatial variables' ${\mathbf y}$ and the `transformed frequency
variables' $\bnu$ by:
$
{\mathbf y}={\mathbf P}_{\bH}\x$ and $\bnu={\mathbf Q}^T_{\bH}\bomega$. Then
$\x={\mathbf Q}_{\bH}\y$
and $\bomega={\mathbf P}_{\bH}^T\bnu$.
For convenience we also define $y=\|\y\|,$ and $\nu=\|\bnu\|,$ to complement
$x=\|\x\|$ (the modulus of the spatial coordinate) and $\omega=\|\bomega\|$.

\begin{figure*}[t]
\centerline{
\includegraphics[scale=0.8]{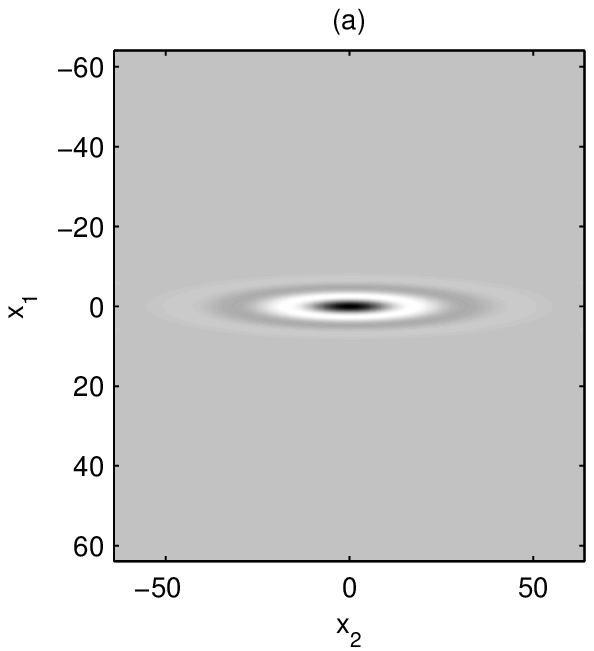}
\includegraphics[scale=0.8]{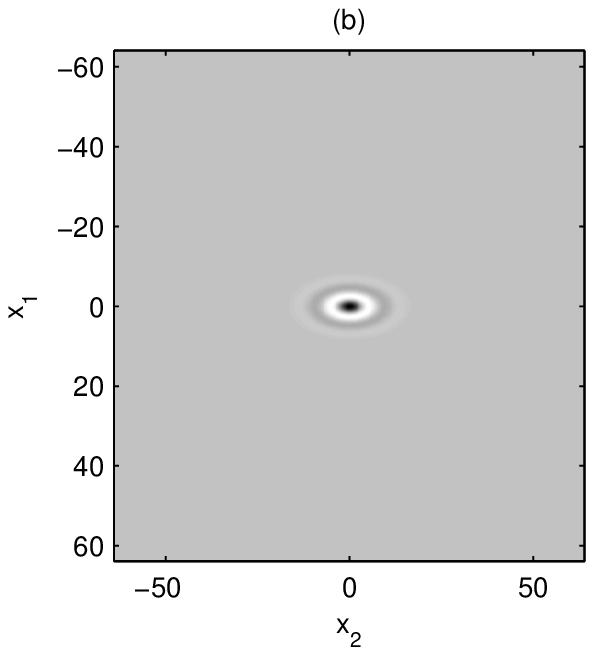}
\includegraphics[scale=0.8]{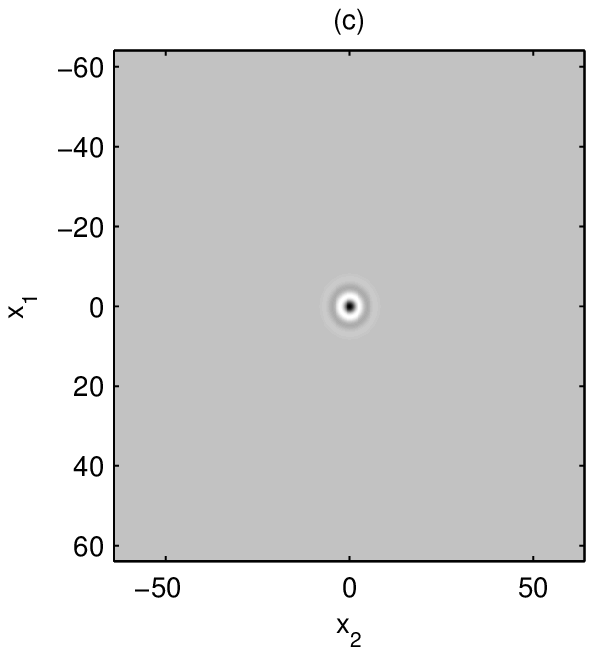}}
\centerline{
\includegraphics[scale=0.8]{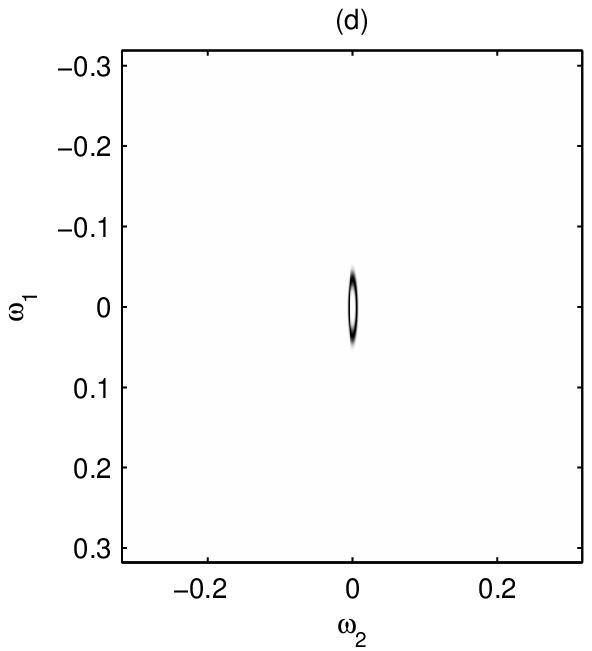}
\includegraphics[scale=0.8]{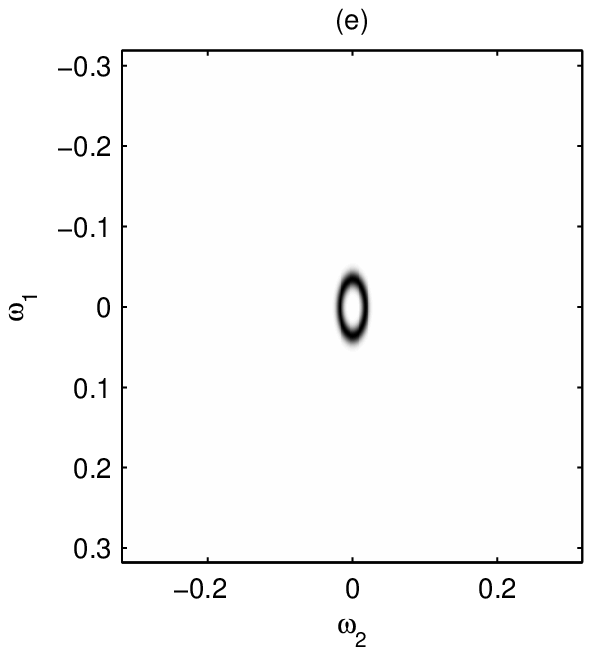}
\includegraphics[scale=0.8]{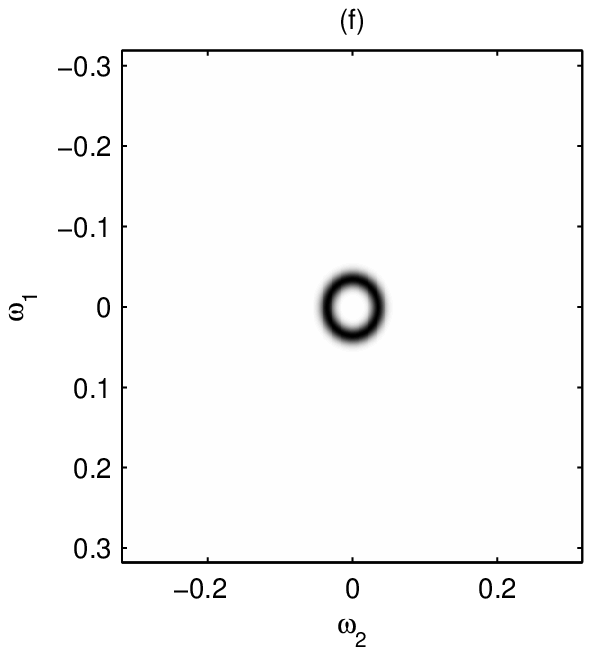}}
\caption{\label{isotropic}
The $n=0$ eigenfunctions in space when
a) $\bP_{\bH}=\left(\left[1 \quad 0\right];\left[0 \quad 0.15\right]\right),$
b) $\bP_{\bH}=\left(\left[1 \quad 0\right];\left[0 \quad 0.5\right]\right)$
and
c) $\bP_{\bH}=\left(\left[1 \quad 0\right];\left[0 \quad 1\right]\right)$,
with $(\beta,\gamma)=(8,3).$ The same functions are also plotted in frequency,
in d), e) and f).
}
\end{figure*}

\subsection{Defining the Geometrically Anisotropic Coherent States \label{defstate}}
We shall now construct a localisation operator that measures the localisation
to regions that spatially decay in $y^2(\bx),$ and decay
in frequency $\bnu(\bomega)$,
away from some set of frequencies.  We first define the fiducial vector, very much in analogue with
choosing a mother wavelet, and from the fiducial vector define a family
of coherent states. 
Define as in Metikas and Olhede \cite{Metikas}[eqn.~4],
$
V^{(1-D,\beta,\gamma)}(\omega)=
\frac{2^{(r+1)/2}\sqrt{\pi\gamma}}{\sqrt{\Gamma(r)}}\omega^{\beta} e^{-\omega^{\gamma}}$
for $\omega>0$, and 
$0$ for $\omega<0$.
We use the radial coherent states used in \cite{Metikas}[eqn.~8]:
\begin{equation}
\label{isomorse}
V^{(\beta,\gamma)}_r(\bomega)=\omega^{-1/2+\gamma/2}V_{1D}^{(1-D,\beta,\gamma)}(\omega),\quad
v^{(\beta,\gamma)}_r(\x)=\int_{-\infty}^{\infty} \int_{-\infty}^{\infty} \frac{\omega^{-1/2+\gamma/2}}{(2\pi)^2}V^{(1-D,\beta,\gamma)}(\omega)
e^{\bj \bomega^T\x}\;d^2\bomega.
\end{equation}
$v^{(\beta,\gamma)}_r(\x)$ is marginally well localised, see 
\cite{Metikas}. 
Define the geometrically anisotropic
fiducial vector by:
\begin{equation}
\label{fiducialvector}
v_{\bH}(\x)=\sqrt{\left|\bP_{\bH}
\right|}v_{r}(\y),\quad V_{\bH}(\bomega)=\left|\bP_{\bH} \right|^{-1/2}
V_r(\mathbf{Q}^T_{\bH}\bomega)=\left|\bP_{\bH} \right|^{-1/2}
V_r(\bnu).\end{equation}
Let $\omega_0=\arg_{\omega>0}\max 
\omega^{-1+\gamma}|V_{1D}^{(1-D,\beta,\gamma)}(\omega)|^2$. $v_{r}^{(\beta,\gamma)}(\x)$ is decaying radially in space
from
$\x={\mathbf{0}},$ and is supported at angular frequencies such that 
$\bm{\omega}=\bm{\omega}_0$, where
$\omega\approx\omega_0$. Thus
$v_{\bH}(\x)$ is local in terms of $\x^T\bH \x$
and
is associated with frequencies $\bomega$ such that 
$\left|\bQ^T_{\bH} \bomega\right|=\omega_0.$

We now need to define the coherent states from the fiducial vector defined
in eqn. (\ref{fiducialvector}).
We define the transformed rotation operator $\tilde{\cal{R}}_{\theta,\bH}$ by:
$\tilde{\cal{R}}_{\theta,\bH} \left\{g\right\}(\x)=g\left(
\tilde{\mathbf{r}}_{-\theta,\bH}\x\right),$ 
 $\tilde{\mathbf{r}}_{-\theta,\bH}=\left[\left({\mathbf{P}}^T_{\bH}
\mathbf{r}_{-\theta}{\mathbf{Q}}^T_{\bH}\right)^T\right]^{-1}.$
Thus if ${\mathbf{P}}_{\bH}$ is itself a rotation, then $
\tilde{\cal R}_{\theta,\bH}$ corresponds to a normal rotation, otherwise $
\tilde{\cal R}_{\theta,\bH}$
defines the act of rotation in the transformed space, i.e. in the $\y$ coordinates.

We define the generalized spatial shift,
also in the transformed space, for any function $g(\x)$ with FT
$G(\bomega)$ by:
\begin{equation}
\label{genspace}
\tilde{\cal T}_{\be,\bH}^{(\gamma)}G(\bomega)=G(\bomega)e^{-\bj \be^T \bnu
\nu^{\gamma-1}},\quad
\tilde{\cal T}_{\be,\bH}^{(\gamma)}g(\x)=\frac{1}{(2\pi)^2}\int_{-\infty}^{\infty}
\int_{-\infty}^{\infty} G(\bomega)e^{\bj\bomega^T(\x-
\nu^{\gamma-1} {\mathbf{Q}}_{\bH}\be)}\;d^2\bomega
.
\end{equation}
Thus if $\gamma=1$ this simplifies to the usual spatial shift, if by $
{\mathbf{Q}}_{\bH}\be$ rather than by $\be$. The generalized spatial shift
is thus adapted to the transformed geometry. If $\gamma
\neq 1$ then
the local effect on $g(\x),$ local to any given wavenumber $\bomega=\bomega_0 ,$ is a shift in space
by $\left|\bP_{\bH} \bomega_0\right|^{\gamma-1}\bQ_{\bH} \be .$
The magnitude of the
shift
depends on $\left|\bP_{\bH} \bomega_0\right|^{\gamma-1}.$ 
The act of a implementing a generalized spatial shift 
corresponds
to shifting the fiducial vector in space by different amounts at different
frequencies, in precise analogue with the radial \cite{Metikas} and 1-D 
\cite{DaubechiesPaul1988} generalized spatial shifts. Introducing the generalized
spatial shift enables us to treat a larger class of localisation regions.
A family of coherent states is then defined by:
\begin{eqnarray}
v_{\bH;\bxi}(\x)&=&
{\cal D}_{a^{1/\gamma}}\tilde{\cal
R}_{\theta,\bH}\tilde{\cal T}_{\be,\bH}^{(\gamma)}v_{\bH}
\left(\x\right),\;
V_{\bH;\bxi}(\bomega)=a^{1/\gamma}
V_{\bH}(a^{1/\gamma}\left(\left[\tilde{\mathbf{r}}_{-\theta,\bH}\right]^{-1}
\right)^T\bomega)e^{-\bj \bomega^T\bQ_{\bH}\be  \nonumber
\|\bQ_{\bH}\bomega\|^{\gamma-1}}\\
&
=& e^{-\bj  \bnu^T\be
\nu^{\gamma-1}}a^{1/\gamma}\left|\bP_{\bH}\right|^{-1/2}V_r(a^{1/\gamma}\rott\bnu)
\label{transformedstate},\;
%
\end{eqnarray}
and $V_{\bH;\bxi}^{(\beta,\gamma)}(\bomega)
=e^{-\left(a\nu+\bj\bnu^T\be\right)\nu^{\gamma-1}}a^{1/\gamma} 2^{\frac{r+1}{2}}\sqrt{\pi \gamma}\left(a^{1/\gamma}\nu\right)^{\frac{2\beta+\gamma-1}{2}}/(
\sqrt{\left|\bP_{\bH}\right|\Gamma(r)}).$

\subsection{The Localisation of the Coherent States \label{locally}}
Before expressing any function in terms of the family of coherent states,
and proving that a resolution of identity may be achieved, we shall determine
the localisation of the coherent states. We need to determine the spatial and spatial frequency locality
of $v_{\bH;\bxi}^{(\beta,\gamma)}(\x)$ as a function of $\bxi,$ to be able
to determine the localisation of the operator. We firstly determine the energy
of the coherent states (see eqn. (\ref{1av}) in the appendix), and noting that
$d^2\bomega=\left|\bP_{\bH}\right| d^2\bnu,$ this then
yields:
\begin{eqnarray}
\nonumber
\langle 1\rangle_{\bH}(\bxi)&=&
\langle v_{\bH;\bxi}^{(\beta,\gamma)},v_{\bH;\bxi}^{(\beta,\gamma)}\rangle
=\frac{1}{(2\pi)^2}\int_{-\infty}^{\infty} 
\int_{-\infty}^{\infty} \left|V_{\bH;\bxi}^{(\beta,\gamma)}
(\bomega)\right|^2\;d^2\bomega 
\nonumber
=\frac{1}{2\pi}\int_{0}^{\infty} 
\omega^{\gamma}\left|V_{1D}^{(\beta,\gamma)2}(\omega)\right|^2\;
d\omega.
\end{eqnarray}
We then calculate the spatial
and spatial frequency average of the function, just like Daubechies and Paul \cite{DaubechiesPaul1988}[eqn.~2.10] . It is easiest to calculate
the transformed coordinate averages, rather than the averages of $\x$ and $\bomega.$ We determine that (see eqn. (\ref{calcfornu}) in the appendix):
\begin{eqnarray}
\nonumber
\langle \nu \rangle_{\bH}(\bxi)&=&\frac{4}{\langle 1\rangle_{\bH}(\bxi)}
\frac{1}{(2\pi)^2}\int_{0}^{\infty} 
\int_{0}^{\infty} \nu \left|V_{\bH;\bxi}(\bomega)\right|^2\;d^2\bomega
\nonumber
=\frac{\left|\bP_{\bH} \right|}{\langle 1\rangle_{\bH}(\bxi)}
\frac{a^{2/\gamma}}{\pi^2}\int_{0}^{\infty} 
\int_{0}^{\infty} \nu \left|V_{\bH}(a^{1/\gamma}\bm{r}_{-\theta}
\bnu)\right|^2\;d^2\bnu\\
&=&\frac{\Gamma\left(r+\frac{1}
{\gamma}+1\right)}{2^{1/\gamma}\Gamma(r+1)}\frac{1}{a^{1/\gamma}}=
\frac{C_3^{(\beta,\gamma)}}{a^{1/\gamma}}
\label{a}
.
\end{eqnarray}
Similarly we can determine that the average spatial position of any member
of the family of coherent states is given (see eqn. (\ref{averagespace22})
in the appendix)
by:
\begin{eqnarray}
\nonumber
\langle y_l \rangle_{\bH}(\bxi)&=&\frac{1}{\langle 1\rangle_{\bH}(\bxi)}
\int_{-\infty}^{\infty} 
\int_{-\infty}^{\infty} y_l \left|v_{\bH;\bxi}(\x)\right|^2\;d^2\x
=  \frac{2^{1/\gamma-2}(\gamma+1)\Gamma(r-1/\gamma+2)}{\Gamma(r+1)}a^{-1+1/\gamma}b_l.
\label{averagespace}
\end{eqnarray}
We define $C_4^{(\beta,\gamma)}=
2^{1/\gamma-2}(\gamma+1)\Gamma(r-1/\gamma+2)/\Gamma(r+1)$, to be the term multiplying $a^{-1+1/\gamma}b_l$.
The average spatial position $\langle \x
\rangle_{\bH}(\bxi)$ has components $
\langle x_l \rangle_{{\bH}}(\bxi)=\left[{\mathbf{Q}}_{\bH}\begin{pmatrix}
\langle y_1 \rangle_{{\bH}}(\bxi)\quad
\langle y_2 \rangle_{{\bH}}(\bxi)
\end{pmatrix} \right]_{l}.$
With this set of relationships we may note the locality of $v_{\bH}^{(\beta,\gamma)}(\x)$ and
use this to construct a suitable localisation region for a function exhibiting
geometric anisotropy. 

\subsection{Defining the Localisation Operator 
\& Resolution of Identity \label{defop}}
\subsubsection{General Operator}
For any $g(\x)\in L^2({\mathbb{R}})$ we define for any region ${\cal A}\subset
{\mathbb{R}}^+\times\left[0,2\pi\right)\times {\mathbb{R}}^2$ the localisation
operator. Define firstly the region of space and spatial frequency given
by:
\begin{equation}
{\cal D}_{\bH}=\left\{(\x,\bomega):\;\;
\x=\langle \x \rangle_{{\bH}}(\bxi),\;\;
\left|\bQ_{\bH}^T\bomega \right|= \langle \nu \rangle_{{\bH}}(\bxi),\quad
\bxi\in {\cal A}
\right\}.
\end{equation}
Then subsequently define the projection operator by:
\begin{eqnarray}
\nonumber
{\cal P}_{{\cal D}_{\bH}}g(\x)&=&
C_{\bH}\int_{{\cal A}} \langle v_{\bH;\bxi},
g\rangle v_{\bH;\bxi}(\x)\;dA_{\bxi}\\
\nonumber
{\cal P}_{{\cal D}_{\bH}}G(\bomega_1)&=&\frac{C_{\bH}}{(2\pi)^2}\int_{{\cal A}} \int_{-\infty}^{\infty}
\int_{-\infty}^{\infty} 
V_{\bH;\bxi}^*(\bomega_2) G(\bomega_2) V_{\bH;\bxi}(\bomega_1)\;d^2
\bomega_2\;dA_{\bxi}\\
\nonumber
{\cal P}_{{\cal D}_{\bH}}^{(\beta,\gamma)}G(\bomega_1)&=&
\frac{C_{\bH}}{(2\pi)^2}\int_{{\cal A}} \int_{-\infty}^{\infty}
\int_{-\infty}^{\infty} 
a^{2/\gamma}V^{(\beta,\gamma)}_r(a^{1/\gamma}\rott\bnu_{1})e^{-\bj \be^T \bnu_{1}
\nu_{1}^{\gamma-1}}\\
&&V^{(\beta,\gamma)*}_r(a^{1/\gamma}\rott\bnu_{2})e^{\bj \be^T \bnu_{2}\nu_{2}^{\gamma-1}}
G\left(\bP^T_{\bH}\bnu_{2}\right)\;d^2\bnu_{2}
\;dA_{\bxi}.
\label{projjy1}
\end{eqnarray}
We now let ${\cal A}\rightarrow {\mathbb{R}}^+\times\left[0,2\pi\right)\times {\mathbb{R}}^2={\cal A}_{\mathrm{all}}$ and
note from Metikas and Olhede \cite{Metikas}[section III.D] that assuming
$
C_{\bH}^{-1}=\frac{1}{2\pi}
\int \int_{\mathbb{R}^2} \nu^{-2\gamma}
\left|V^{(\beta,\gamma)}_r(\bnu) \right|^2\;d^2\bnu < \infty,
$
we have that
$
{\cal P}_{{\cal D}_{\mathrm{all}}}^{(\beta,\gamma)}G(\bomega_1)
=G\left(\bP^T_{\bH}\bnu_{1}\right)=G\left(\bomega_1\right).
$
We thus achieve a `resolution of identity'.
This demonstrates that eqn. (\ref{projjy1}) is a suitable localisation operation. If we include all space when integrating then we retrieve the full function
$g(\x)$. Of course, to appropriately define a localised function, the operator needs to be specified very carefully.
 
We start by defining the appropriate restriction of the local index $\bm{\xi}$
for a fixed $C>1$
by:
\begin{equation}
\label{geoanisoregion}
{\cal A}(C)=
\left\{(a,\theta,\be):\quad a^2+b^2+1\le 2a C,\quad \theta\in\left[0,2\pi\right),\quad
b=\sqrt{b_1^2+b_2^2}
\right\}.
\end{equation}
Note that
${\cal A}(C)$ does not depend on $\bH,$ the transformation of the distance metric,
or the shape parameters $(\beta,\gamma).$
We shall now use the coherent states
defined by eqn. (\ref{transformedstate}) so that ${\cal A}(C)$ of equation (\ref{geoanisoregion}) gets mapped
onto a region of space and spatial frequency  that is given by:
\begin{equation}
\label{geoaniso1}
{\cal D}_{\bH}^{(\beta,\gamma)}(C)=\left\{(\x,\bomega):\quad \left(y \nu E_3^{(\beta,\gamma)}\right)^2+\left(
\nu^{\gamma}-CE_4^{(\beta,\gamma)}\right)^2 \le E_4^{(\beta,\gamma)2}\left(C^2-1 \right)
\right\},
\end{equation} 
with $C_3^{(\beta,\gamma)},$
and $C_4^{(\beta,\gamma)},$
constants whose values are given by eqns. (\ref{a}) and (\ref{averagespace}).
For simplicity we have also defined:
\begin{equation}
E_{3}^{(\beta,\gamma)}=\left(\frac{C_3^{(\beta,\gamma)(\gamma-1)}}{C_4^{(\beta,\gamma)}}\right),\quad
E_4^{(\beta,\gamma)}=C_3^{(\beta,\gamma)\gamma}.
\end{equation}
Thus the region ${\cal D}_{\bH}^{(\beta,\gamma)}(C)$ is specified by four
different constituent parts: the transformation $\bP_{\bH},$ the shape parameters
$(\beta,\gamma)$ and $C$, the hypervolume parameter.

We note from eqn. (\ref{geoanisoregion}) that
$C-\sqrt{C^2-1}\le a \le C+\sqrt{C^2-1}$ and so 
the maximum and minimum wave number in ${\cal D}_{\bH}^{(\beta,\gamma)}(C)$
are therefore given by:
$
\nu_{\max}=C_3^{(\beta,\gamma)}/\left[\sqrt[\gamma]{\left(C-\sqrt{C^2-1}\right)}\right]$,
$\nu_{\min}=C_3^{(\beta,\gamma)}/\left[\sqrt[\gamma]{\left(C+\sqrt{C^2-1}\right)}\right]$.
The local transformation $\bP_{\bH}$ specifies the local
region ${\cal D}_{\bH}^{(\beta,\gamma)}(C)$ in terms of the anisotropy between
$x_1$ and $x_2$, as well as $\omega_1$ and $\omega_2$,
as the
space is parameterised in $\y$ and $\bnu.$
The shape parameters $(\beta,\gamma)$ determine the compromise between $\y$ and $\bnu,$
and $C$ determines the hypervolume of the region.
We reparameterise
the region ${\cal D}_{\bH}^{(\beta,\gamma)}(C)$ in  $\y$ and $\bnu,$ as ${\cal D}_{\bH,y}^{(\beta,\gamma)}(C).$
For a fixed value of $C,$ the hypervolume of 
${\cal D}_{\bH}^{(\beta,\gamma)}(C),$ is given by:
\begin{eqnarray}
\nonumber
\left|{\cal D}^{(\beta,\gamma)}_{\bH}(C)\right|&=&
\int_{{\cal D}^{(\beta,\gamma)}_{\bH}(C)}\;d^2\x\;d^2\bomega=
\int_{{\cal D}^{(\beta,\gamma)}_{\bH,y}(C)}\;d^2\y\;d^2\bnu
=4 \pi^2 \int_{{\cal D}_{\bH,y}^{(\beta,\gamma)}(C)} y\;dy\;\nu\;d\nu.
\nonumber
\end{eqnarray}
We change variables to
$
a=\frac{C_3^{(\beta,\gamma)\gamma}}{\nu^{\gamma}},$
$b_l=\frac{C_3^{(\beta,\gamma)(\gamma-1)}y_l}{C_4^{(\beta,\gamma)}\nu^{\gamma-1}},$
and $b=\sqrt{b_1^2+b_2^2}=\frac{C_3^{(\beta,\gamma)(\gamma-1)}y}{C_4^{(\beta,\gamma)}\nu^{\gamma-1}}.$
Then:
\begin{eqnarray}
\nonumber
A(\beta,\gamma){\mathrm{Area}}(C)&=&\left|{\cal D}^{(\beta,\gamma)}_{\bH}(C)\right|=
4 \pi^2 \int_{C-\sqrt{C^2-1}}^{C+\sqrt{C^2-1}}\int_{0}^{2aC-1-a^2} 
\frac{C_4^{(\beta,\gamma)}b}{a^{1-1/\gamma}}
\frac{C_3^{(\beta,\gamma)}}{a^{1/\gamma}}\frac{C_3^{(\beta,\gamma)}
C_4^{(\beta,\gamma)}}{\gamma}\;\frac{da\;db}{a^2}\\
&=&\frac{(\gamma+1)^2\Gamma^2(r-\frac{1}{\gamma}+2)
\Gamma^2(r+\frac{1}{\gamma}+1)}{2^{5}\Gamma^4(r+1)\gamma}
\left[2C\sqrt{C^2-1}
+\log\left(\frac{C-\sqrt{C^2-1}}{C+\sqrt{C^2-1}}\right)
\right] \nonumber .
\end{eqnarray}
We obtain the form starting from eqn. (\ref{areaappendix}) in the appendix.
This defines $ A(\beta,\gamma)=\frac{(\gamma+1)^2\Gamma^2(r-\frac{1}{\gamma}+2)
\Gamma^2(r+\frac{1}{\gamma}+1)}{2^{5}\Gamma^4(r+1)\gamma}$ and
${\mathrm{Area}}(C)=\left[2C\sqrt{C^2-1}
+\log\left(\frac{C-\sqrt{C^2-1}}{C+\sqrt{C^2-1}}\right)
\right].$ Note that $A(\beta,\gamma)$ is only a function of $(\beta,\gamma)$
via $r$
and so does not depend on $C,$ whilst ${\mathrm{Area}}(C)$ is only a function
of $C$ and does not change with $(\beta,\gamma).$ Also $\frac{d}{dC}
{\mathrm{Area}}(C)>0$ and so there is a 1-1 map between C and the area.

Metikas and Olhede \cite{Metikas}[eqn.~22]
defined a genuinely 2-D localisation operator, but were only able to determine
its approximate eigenfunctions. To be able to exactly derive the eigenfunctions,
they defined an operator only valid for radial functions, in  
\cite{Metikas}[eqn.~12]. Similarly, we intend
to define an operator only valid for geometrically anisotropic functions.
We assume that $g(\x)$ is a geometrically anisotropic function, and write
it as $g_A(\x)$.
Then in analogue with \cite{Metikas} we
define the geometrically anisotropic coherent state in terms of the simpler
indexing of
$\bxi_I=\left(a,b\right)^T$
by:
\begin{eqnarray}
V_{\bH,\bxi_I}^{(1,\beta,\gamma)}(\bomega)&=&
\frac{a^{1/\gamma} \left(a^{1/\gamma} \nu\right)^{-1/2+\gamma/2}
V^{(1-D,\beta,\gamma)}(a^{1/\gamma} \nu) 
\cos\left(\nu^{\gamma}b-\pi/4 \right)}{\sqrt{2^{-1}\pi\nu^{\gamma}b
\left|\bP_{\bH}\right|}},\\
v_{\bH,\bxi_I}^{(1,\beta,\gamma)}(\x)&=&
\frac{1}{(2\pi)^2}\int \int_{{\mathbb{R}}^2}
V_{\bH,\bxi_I}^{(1,\beta,\gamma)}(\bomega)\;e^{i\bomega^T\x}\;d^2\bomega
.
\end{eqnarray}
As we may consider $v_{\bH,\bxi_H}^{(1,\beta,\gamma)}(\x)$ as the coherent
state $v_{\bH,\bxi}^{(\beta,\gamma)}(\x)$ that has been averaged across directions in space over $b$ constant (see \cite{Metikas}), we may consider
$v_{\bH,\bxi_H}^{(1,\beta,\gamma)}(\x)$ local to the set of spatial positions
\begin{eqnarray}
\x=\bQ_{\bH}\y:\quad y^2=\left|\langle \by\rangle_{\bH} (\bxi)\right|^2,\quad
\bomega=\bQ_{\bH}\bnu:\quad \nu=\langle \nu\rangle_{\bH} (\bxi),
\end{eqnarray}
noting that $\left|\langle \by\rangle_{\bH} (\bxi)\right|^2$ only depends on
$a$ and $b$ (see eqn. (\ref{averagespace})), whilst $\langle \nu\rangle_{\bH} (\bxi)$ only depends on $a$ (see eqn. (\ref{a})).
By only using $v_{\bH,\bxi_I}^{(1,\beta,\gamma)}(\x) $ for
$\bxi_H\in {\cal A}_r(C)$
with
$
{\cal A}_r(C)=
\left\{(a,\theta,\be):\; a^2+b^2+1\le 2a C
\right\}$
we construct a function local in space and spatial frequency to
\begin{equation}
{\cal{D}}_{r,\bH}(C)=\left\{(\x,\bomega):\;\;
\left|\bP\x\right|=\left|\langle \y \rangle_{{\bH}}(\bxi)\right|,\;\;
\left|\bQ_{\bH}^T\bomega \right|= \langle \nu \rangle_{{\bH}}(\bxi),\quad
b_1^2+b_2^2=b^2,\quad
(a,b)\in {\cal A}_r(C)
\right\}.
\nonumber
\end{equation}

In analogue with Metikas and Olhede \cite{Metikas} we then may define the localisation operator with
$\langle g_1,g_2 \rangle_{\bH}=(2\pi)^{-1}\int G_1^*(\omega)G_2(\omega)\;\nu\;
\left|\bP_{\bH}\right|\;d\nu,$
by:
\begin{equation}
{\cal P}_{{\cal{D}}_{r,\bH}(C)}^{(\beta,\gamma)}\left\{G_A\right\}(\bomega)=
C_{\bH,A}\int \int_{{\cal{A}}_r(C)}
V_{\bH,\bxi_H}^{(1,\beta,\gamma)}(\bomega)\langle
v_{\bH,\bxi_H}^{(1,\beta,\gamma)}, g_A
\rangle_{\bH} \frac{da}{a^3}b\;db.
\label{anisoaver}
\end{equation}
Clearly by appropriate change of variables, we can note directly from Metikas
and Olhede \cite{Metikas}[eqn.~14] that a resolution of identity may be achieved
as long as 
$C_{\bH,A}=(r-1)/2.$ 

\section{Determining the Eigensystem \label{eigen}}
We now intend to demonstrate that the eigenfunctions of the proposed operators
can in the instance of eqn. (\ref{anisoaver}) be determined exactly, or in the
instance of eqn. (\ref{projjy1}) be determined approximately.

Firstly 
we note that by eqn. (\ref{projjy1}) that determining the eigenfunctions of the operator defined from the Morse coherent states for an arbitrary function
corresponds to for a given $\bP_{\bH},$ $(\beta,\gamma)$ and $C$ solving the set of eqns. $
{\cal P}_{{\cal D}_{\bH}^{(\beta,\gamma)}(C)}\Psi(\bomega_1)=\lambda \Psi(\bomega_1),
$ or:
\begin{eqnarray}
\label{eigenequation}
{\cal P}_{{\cal D}_{\bH}(C)}^{(\beta,\gamma)}\Psi(\bomega_1)&=&
\frac{C_{\bH}}{(2\pi)^2}\int_{{\cal A}(C)} \int_{-\infty}^{\infty}
\int_{-\infty}^{\infty} a^{2/\gamma}
V^{(\beta,\gamma)}_r(a^{1/\gamma}\rott\bnu_1)e^{-\bj \be^T \bnu_1
\nu_1^{\gamma-1}}\\
\nonumber
&&V^{(\beta,\gamma)*}_r(a^{1/\gamma}\rott\bnu_2)e^{\bj \be^T \bnu_2\nu_2^{\gamma-1}}
\Psi(\bomega_{2})\;\frac{d^2\bomega_{2}}{\left|\bP_{\bH}\right|}
\;dA_{\bxi}\\
\nonumber
&=&\frac{C_{\bH}}{(2\pi)^2}\int_{{\cal A}(C)} \int_{-\infty}^{\infty}
\int_{-\infty}^{\infty} a^{2/\gamma}V^{(\beta,\gamma)}_r
(a^{1/\gamma}\rott{\mathbf{Q}}^T_{\bH}
\bomega_{1})e^{-\bj \be^T {\mathbf{Q}}^T_{\bH}\bomega_{1}
\left|{\mathbf{Q}}^T_{\bH}\bomega_{1}\right|^{\gamma-1}}\\
&&V^{(\beta,\gamma)*}_r(a^{1/\gamma}\rott\bnu_2)e^{\bj \be^T \bnu_2\nu_2^{\gamma-1}}
\Psi({\mathbf{P}}^T_{\bH}\bnu_2)\;d^2\bnu_2
\;dA_{\bxi}.
\nonumber
\end{eqnarray}
Of course we note that $\bomega_1=\bP^T_{\bH}\bnu_1,$ and so we find that:
\begin{eqnarray}
\nonumber
{\cal P}_{{\cal D}_{\bH}(C)}^{(\beta,\gamma)}\Psi({\mathbf{P}}^T_{\bH}\bnu_1)
&=&\frac{C_{\bH}}{(2\pi)^2}\int_{{\cal A}(C)} \int_{-\infty}^{\infty}
\int_{-\infty}^{\infty} a^{2/\gamma}V^{(\beta,\gamma)}_r
(a^{1/\gamma}\rott
\bnu_{1})e^{-\bj \be^T \bnu_{1}
\left|\bnu_{1}\right|^{\gamma-1}}\\
&&V^{(\beta,\gamma)*}_r(a^{1/\gamma}\rott\bnu_2)e^{\bj \be^T \bnu_2\nu_2^{\gamma-1}}
\Psi({\mathbf{P}}^T_{\bH}\bnu_2)\;d^2\bnu_2
\;dA_{\bxi}.
\label{eigen11}
\end{eqnarray}
Let $l=\beta-\frac{1}{2}$ and $m=\gamma$.
Comparing eqn. (\ref{eigen11}) with \cite{Metikas}[eqns.~7,~21],
and denoting the eigenfunctions of the operator defined 
by \cite{Metikas}[eqns.~7,~21] as
$\Psi_{n;l,m}^{(e)}(\bomega)$,
we obtain that
the eigenfunctions of equation (\ref{eigen11}), denoted by $\psi_{n;l,m,\bH}^{(A)}(\bomega)$, are given by the same functional
form,
but where the argument has been adjusted to the geometric anisotropy.
$\psi_{n;l,m,\bH}^{(A)}(\bomega)$ have the same eigenvalues as the isotropic Morse wavelets,
and we denote them $\left\{\lambda_{n;r,\bH}(C)\right\}$.
We may then note that the geometrically anisotropic eigensystem is given
by:
\begin{equation}
\label{eigensystem}
\Psi_{n;l,m,\bH}^{(A)}(\bomega)=\Psi_{n;l,m}^{(e)}(\bQ_{\bH}^T \bomega),\quad
\lambda_{n;r,\bH}(C)=\lambda_{n,r}(C).
\end{equation}
Thus, whenever a set of eigenfunctions are determined for the operator of
eqns. (7) and (21) in \cite{Metikas}, these automatically correspond
to eigenfunctions of the geometrically anisotropic localisation operator, once
the argument has been adjusted.
%
Determining the exact eigenfunctions of eqn. (\ref{eigen11}) is in general
not an analytically tractable problem. In analogue with Metikas and Olhede
\cite{Metikas} we instead determine the eigenfunctions of the geometrically anisotropic localisation operator:
\begin{eqnarray}
\label{aniso1}
{\cal P}_{{\cal{D}}_{r,\bH}(C)}^{(\beta,\gamma)}\left\{\psi_A\right\}(\x)
&=&\lambda\psi_A(\x)
\\
\nonumber
{\cal P}_{{\cal{D}}_{r,\bH}(C)}^{(\beta,\gamma)}\left\{\Psi_A\right\}(\bomega_1)&=&
C_{\bH,A}\int \int_{{\cal{A}}_r(C)}
V_{\bH,\bxi_I}^{(1,\beta,\gamma)}(\bomega)\langle
v_{\bH,\bxi_I}^{(1,\beta,\gamma)}, g_A
\rangle_{\bH} \frac{da}{a^3}b\;db\\
\nonumber
&=& C_{\bH,A}\int_0^{\infty} \int \int_{{\cal{A}}_r(C)} a^{1/\gamma}(\nu_1\nu_2)^{-1/2+\gamma/2}
V^{(1-D,\beta,\gamma)}(a^{1/\gamma} \nu_1) 
\frac{\cos\left(\nu^{\gamma}_1b-\pi/4 \right)}{\sqrt{2^{-1}\pi\nu_1^{\gamma}b}}
\\
&&\Psi_A\left(\bP_{\bH}^T\bnu_2\right)
V^{(1-D,\beta,\gamma)}(a^{1/\gamma} \nu_2) 
\frac{\cos\left(\nu_2^{\gamma}b-\pi/4 \right)}{\sqrt{2^{-1}\pi\nu_2^{\gamma}b}}\;\nu_2\;d\nu_2
\frac{da\; b\;db}{a^2}.
\label{anisoaver2}
\end{eqnarray}
Clearly comparing this with eqn. (13) of \cite{Metikas}, we determine
that the eigensystem is given by:
$
\Psi_{n;l,m,\bH}^{(A)}(\bomega)=\Psi_{n;l,m}^{(e)}(\bQ_{\bH}^T \bomega)$
and
$\lambda_{n;r,\bH}(C)=\lambda_{n,r}(C),$
where we note from Metikas and Olhede \cite{Metikas} that the eigenfunctions
and eigenvalues take the form:
\begin{eqnarray*}
\Psi_{n;l,m}^{(A)}(\bomega)=\sqrt{2}A_{n;l+1/2,m}\nu^l e^{-\nu^m}
L_n^{c}(2\nu^m),\quad
\lambda_{n,r}(C)=\frac{\Gamma(r+n)}{\Gamma(N+1)\Gamma(r-1)}\int_0^{\frac{C-1}{C+1}}
x^n \left(1-x\right)^{r-2}\;dx,
\end{eqnarray*}
where $c=(2l+2)/m-1$ and $L_n^c(\cdot)$ denotes the generalized Laguerre polynomial
\cite{AS}[p.~783]. 

As an example of geometrical anisotropy
for $0<\varepsilon \ll 1$ we take:
$\bP_{\bH}(\varepsilon)=\left[1 \quad 0;\;
0 \quad \varepsilon
\right],$ $\bQ_{\bH}(\varepsilon)=\left[1 \quad 0;\;
0 \quad \varepsilon^{-1}
\right],$ $ \bH(\varepsilon)=
\left[
1 \quad 0;\;
0 \quad \varepsilon^2  \right]$.
We note that $
\omega=\sqrt{\omega_1^2+\varepsilon^{-2}\omega_2^2}.$
We obtain increasing anisotropy as $\varepsilon $ approaches zero.  

In the isotropic limit we obtain the results of Metikas and Olhede \cite{Metikas}
and with
$\bP_{\bH}=\eye,$
the isotropic Morse wavelets are retrieved \cite{Metikas}. Compare the range of possible
localisation regions that may be found in this class: see Figure \ref{isotropic}. As we do not change the values of $\beta$ and $\gamma$
these functions have the same localisation value for a fixed value of 
$C$. As we are for the highly anisotropic case extending the extent of the
function in one direction, trivially, to conserve the area size, it is compressed
in the corresponding canonical variable, as is apparent from the frequency domain plot.
As the spatial extent extends in $x_2$, the function squashes towards $f_2=0$.
The proposed framework allows the quantification of
the stretching explicitly, and the discussion of both $\varepsilon\rightarrow
0$ and $\varepsilon\rightarrow \infty$.
 
\section{Conclusions \label{con}}
This correspondence has defined geometric anisotropic localisation, and derived
the optimally local functions for geometrically
anisotropic and nonseparable regions of space and spatial frequency.
A full set of eigenfunctions with associated eigenvalues
have been determined for this class, useful for the study of non-stationary
fields in geostatistics. 


\appendix
Renormalising the functions
to unit energy using eqn.
(\ref{transformedstate}) we have:
\begin{eqnarray}
\label{1av}
\langle 1 \rangle_{\bH}\left(\bxi\right)&=&\langle 
v_{\bH;\bxi}^{(\beta,\gamma)},v_{\bH;\bxi}^{(\beta,\gamma)}
\rangle 
=
(2\pi)^{-2}\int_{-\infty}^{\infty} 
\int_{-\infty}^{\infty} V_{\bH,\bxi}^{(\beta,\gamma)*}\left(\bomega\right)
V_{\bH,\bxi}^{(\beta,\gamma)}\left(\bomega\right)\;d^2\bomega
\\
\nonumber
&=&(2\pi)^{-2}\int_{-\infty}^{\infty} 
\int_{-\infty}^{\infty} 
\frac{\left|a^{1/\gamma}V^{(\beta,\gamma)}_r(a^{1/\gamma}
\rott\bnu)\right|^2\;d^2\bomega}{\left|\bP_{\bH}\right|}
%
=(2\pi)^{-2}\int_{-\infty}^{\infty} 
\int_{-\infty}^{\infty}\omega^{-1+\gamma}\left|V_{1D}^{(\beta,\gamma)2}(\omega)\right|^2\;
d^2\bomega.
\end{eqnarray}
\begin{eqnarray}
\nonumber
\langle \nu \rangle_{\bH}(\bm{\xi})
&=&\frac{1}{\langle 1\rangle_{\bH}(\bm{\xi})}
\frac{1}{(2\pi)^2 a^{1/\gamma}}\int_{-\infty}^{\infty} 
\int_{-\infty}^{\infty} \nu \left|V_{\bH}^{(\beta,\gamma)}(\bm{\nu})\right|^2\;d^2\bm{\nu}
=\frac
{2}{r}
\frac{2^{r+1}\pi \gamma}{\Gamma(r)}
\frac{1}{2\pi a^{1/\gamma}}\int_{0}^{\infty} 
\omega^{\gamma+1} \omega^{2\beta} e^{-2\omega^{\gamma}}\;d\omega\nonumber
\\
&=&\frac{1}{\Gamma(r+1)}\frac{2^{r+1}}{a^{1/\gamma}}\int_{0}^{\infty}
s^{(2\beta+2)/\gamma} e^{-2s}\;ds=\frac{\Gamma\left(r+\frac{1}
{\gamma}+1\right)}{2^{r+\frac{1}
{\gamma}+1}\Gamma(r+1)}\frac{2^{r+1}}{a^{1/\gamma}}=\frac{C_3^{(\beta,\gamma)}}{a^{1/\gamma}}.
\label{calcfornu}\\
\nonumber
\langle y_l \rangle_{\bH}(\bm{\xi})&=&\frac{1}{\langle 1\rangle_{\bH}(\bm{\xi})}
\int_{-\infty}^{\infty} 
\int_{-\infty}^{\infty} y_l \left|v_{\bH;\bm{\xi}}^{(\beta,\gamma)}(\x)\right|^2\;d^2\x
=\bj
\int_{-\infty}^{\infty} 
\int_{-\infty}^{\infty} \left[\frac{\partial}{\partial \nu_l} 
V_{\bH;\bm{\xi}}^{(\beta,\gamma)}(\bomega)\right]
\frac{V_{\bH;\bm{\xi}}^{(\beta,\gamma)*}(\bomega)}{(2\pi)^2\langle 1\rangle_{\bH}(\bm{\xi})}
\;d^2\bomega\nonumber\\
&=&\frac{2^r \gamma a^{r+1}}{\Gamma(r+1)\pi}\int_{-\infty}^{\infty}
\int_{-\infty}^{\infty} 
\left[b_l \omega^{\gamma-1}+ (\gamma-1)b_l\omega^{\gamma-2}
\frac{\omega_l^2}{\omega}
\right]\omega^{2\beta+(\gamma-1)}e^{-2a\omega^{\gamma}}
\;d^2\bm{\omega}\nonumber\\
&=&\frac{ 2^{r}  a^{r+1} b_l(\gamma+1)}{\Gamma(r+1)} \int_{0}^{\infty}
u^{2\beta/\gamma+2-1/\gamma}e^{-2au }
\;du u^{1/\gamma-1}\nonumber
=  \frac{2^{1/\gamma-2}(\gamma+1)\Gamma(r-1/\gamma+2)a^{-1+1/\gamma}b_l}{\Gamma(r+1)}.
\label{averagespace22}
\end{eqnarray}
We define $C_4^{(\beta,\gamma)}$ by $C_4^{(\beta,\gamma)}=
\langle y_l \rangle_{\bH}(\bm{\xi})/(b_l a^{1/\gamma-1})$
Furthermore the localisation region has size:
\begin{eqnarray}
\nonumber
\left|{\cal D}_{\bH}^{(\beta,\gamma)}(C)\right|&=&
4 \pi^2 \int_{C-\sqrt{C^2-1}}^{C+\sqrt{C^2-1}}\int_{0}^{2aC-1-a^2} 
\frac{C_4^{(\beta,\gamma)}b}{a^{1-1/\gamma}}
\frac{C_3^{(\beta,\gamma)}}{a^{1/\gamma}}\frac{C_3^{(\beta,\gamma)}
C_4^{(\beta,\gamma)}}{\gamma}\;\frac{da\;db}{a^2}\\
&=&\frac{2 \pi^2\left(C_3^{(\beta,\gamma)}C_4^{(\beta,\gamma)} \right)^2}{\gamma}
\int_{C-\sqrt{C^2-1}}^{C+\sqrt{C^2-1}}\left(2aC-1-a^2\right)\frac{da}{a^3}
\nonumber\\ 
&=&\frac{(\gamma+1)^2\Gamma^2(r-\frac{1}{\gamma}+2)
\Gamma^2(r+\frac{1}{\gamma}+1)}{2^{5}\Gamma^4(r+1)\gamma}
\left[2C\sqrt{C^2-1}
+\log\left(\frac{C-\sqrt{C^2-1}}{C+\sqrt{C^2-1}}\right)
\right] \nonumber \\
&=&A(\beta,\gamma){\mathrm{Area}}(C).
\label{areaappendix}
\end{eqnarray}

\end{document}